\documentclass{article}

\usepackage{arxiv}

\usepackage[utf8]{inputenc} 
\usepackage[T1]{fontenc}    
\usepackage{hyperref}       
\usepackage{url}            
\usepackage{booktabs}       
\usepackage{amsfonts}       
\usepackage{nicefrac}       
\usepackage{microtype}      
\usepackage{graphicx}
\usepackage{natbib}
\usepackage{doi}
\usepackage{subcaption}
\usepackage{caption}

\usepackage{amssymb,amsthm,amsmath}
\usepackage{xcolor,paralist,hyperref,titlesec,fancyhdr,etoolbox}
\newtheorem{theorem}{Theorem}[section]
\newtheorem{definition}{Definition}[section]
\newtheorem{example}{Example}[section]
\newtheorem{lemma}{Lemma}[section]
\newtheorem{proposition}{Proposition}[section]
\newtheorem{corollary}{Corollary}[section]

\title{Universal Approximation of Parametric Optimization via Neural Networks with Piecewise Linear Policy Approximation}

\author{Hyunglip Bae \\
	Department of Industrial and System Engineering\\
	KAIST\\
	South Korea \\
	\texttt{qogudflq@kaist.ac.kr} \\
	\And
	Jang Ho Kim\\
	Depart of Industrial and Management System Engineering\\
    Department of Big Data Analytics\\
	Kyung Hee University\\
	South Korea\\
	\texttt{janghokim@khu.ac.kr} \\
    \And
    Woo Chang Kim\thanks{Corresponding Author}\\
    Department of Industrial and System Engineering\\
    KAIST\\
    South Korea\\
    \texttt{wkim@kaist.ac.kr}
}

\renewcommand{\shorttitle}{Universal Approximation of Parametric Optimization}

\hypersetup{
pdftitle={A template for the arxiv style},
pdfsubject={q-bio.NC, q-bio.QM},
pdfauthor={David S.~Hippocampus, Elias D.~Striatum},
pdfkeywords={First keyword, Second keyword, More},
}

\begin{document}
\maketitle

\begin{abstract}
Parametric optimization solves a family of optimization problems as a function of parameters.
It is a critical component in situations where optimal decision making is repeatedly performed for updated parameter values, but computation becomes challenging when complex problems need to be solved in real-time. 
Therefore, in this study, we present theoretical foundations on approximating optimal policy of parametric optimization problem through Neural Networks and derive conditions that allow the Universal Approximation Theorem to be applied to parametric optimization problems by constructing piecewise linear policy approximation explicitly. 
This study fills the gap on formally analyzing the constructed piecewise linear approximation in terms of feasibility and optimality and show that Neural Networks (with ReLU activations) can be valid approximator for this approximation in terms of generalization and approximation error.
Furthermore, based on theoretical results, we propose a strategy to improve feasibility of approximated solution and discuss training with suboptimal solutions.
\end{abstract}

\keywords{Neural Networks \and Universal Approximation \and Parametric Optimization \and The Maximum Theorem}

\section{Introduction}\label{intro}
Consider a parametric optimization problem  parameterized by $\theta$
\[
\min_{x}\quad 
f(x,\theta) \quad \text{subject\;to} \quad x\in C(\theta).
\]
where $x$ is decision variable, $f$ is objective function and $C(\theta)$ is feasible region for $x$. 
Parametric optimization involves a process of solving a family of optimization problems as a function of parameters. (\citet{nikbakht2020unsupervised, still2018lectures}).
Therefore, it is commonly applied when decisions are made repeatedly as the parameters change, while the fundamental problem structure remains constant over the entire duration.
The parameters are often determined by the environment where decision makers typically can not have control. 
Therefore, an optimization problem must be solved after observing the current state of the environment over and over. 
As solving an optimization problem requires certain computational time, it inevitably causes delays between repetitive decisions, especially for large-scale and complex optimization problems.

There are many application fields that parametric optimization plays a significant role, including robotics (\citet{khalaf2016parametric}), autonomous vehicle control (\citet{daryina2021parametric}), supply chain optimization (\citet{bai2016robust}), and energy system management (\citet{wang2014parametric}).
For example, an autonomous vehicle needs optimal decisions that depend on changes in speed, road conditions, or amount of traffic. 
Any delay in decision-making for autonomous vehicle control could lead to mishandle or even serious traffic accidents.
Similarly, various decisions based on optimization are required for managing responsive and adaptive manufacturing systems. 
Sequential and interconnected systems magnify the importance of computation speed as well as minimal error. 
Real-time decision making is also crucial in high frequency trading of financial assets. 
Delays in trade execution, even for a fraction of a second, can lead to significant losses for financial management firms. 
These applications clearly highlight the importance of latency issues in parametric optimization.

In situations where a family of optimization problems needs to be solved repeatedly, the following two characteristics can be observed.
First, the structure of optimization problems that are solved repeatedly is identical except for input parameters, which means the dependent variable for the optimal policy is input parameters. Second, input parameters and their corresponding optimal policy are accumulated as optimization problems are solved for new input parameters. 
This second case gives potential for supervised learning. 
Thus, it is intuitive and beneficial to approximate the mapping from input parameters to the optimal policy via Machine Learning (ML) techniques in that it is efficient and scalable.

Therefore, in this study, we focus on applying Neural Networks (NN) to universally approximate parametric optimization problems.
We build theoretical foundations on approximating direct mapping from input parameters to optimal solution through NN universally and derive conditions that allow the Universal Approximation Theorem (UAT) to be applied to parametric optimization problems by constructing piecewise linear policy approximation explicitly.
More specifically, we cast single-valued continuous piecewise linear approximation for optimal solution of parametric optimization and analyze it in terms of feasibility and optimality and show that NN with ReLU activations can be valid approximator in terms of generalization and approximation error.
There are various works on the expressive power of NN for approximating functions, however,  to the best of our knowledge,  existing literature lacks theoretical analysis on the applicability of UAT when the target function of NN is the result from parametric optimization problem, and our study is the first to fill this gap.

\subsection{Related Work.}\label{relwork}
The power of NN as a universal approximator has been extensively validated over several decades. 
Pointedly, initial work on the UAT show that, for any continuous function on a compact set, there exist a feedforward NN with a single hidden layer the uniformly approximates the function arbitrarily well (\citet{hornik1989multilayer, cybenko1989approximation, funahashi1989approximate, barron1993universal}).
Recently, there has been a growing interest in exploring the capability of NN for approximating functions stemmed from the works by \citet{liang2016deep} and \citet{yarotsky2017error};  for more recent developments, see also \citet{yarotsky2018optimal, petersen2018optimal, shaham2018provable, shen2019deep, daubechies2022nonlinear, lu2021deep}.
From theoretical view, \citet{telgarsky2016benefits} discussed the benefits of depth in NN, which led to various research on arbitrary depth (\citet{lu2017expressive, hanin2017approximating, kidger2020universal, Sun_Chen_Wang_Liu_Liu_2016, daniely2017depth}). 
There are also numerous extensions of the UAT that is derived from other networks (\citet{baader2019universal, lin2018resnet}) or is generalized to unbounded activation functions (\citet{sonoda2017neural}), discontinuous activation functions (\citet{leshno1993multilayer}), non-compact domains (\citet{kidger2020universal}), interval approximation (\citet{wang2022interval}), distribution approximation (\citet{lu2020universal}) and invariant map (\citet{yarotsky2022universal}). 

Stability analysis of optimization problems plays an important role in control theory which has been utilized in many applications such as electronic engineering (\citet{wang2019parametric}), biology (\citet{motee2012stability}), and computer science (\citet{bubnicki2005modern}).
\citet{berge1963espaces} first proved the Maximum Theorem that provides conditions for the continuity of optimal value function and upper hemicontinuity of optimal policy with respect to its parameters. 
Since the Maximum Theorem only guarantees upper hemicontinuity of optimal policy, this led to extensions on studying conditions for lower hemicontinuous of optimal policy. 
Approaches in such literature are largely divided into two main types. 
Some provide conditions for lower hemicontinuous directly (\citet{robinson1974sufficient}, \citet{zhao1997lower}, \citet{kien2005lower}). 
On the other hand, others provide conditions by limiting structure to be linear (\citet{bohm1975continuity}, \citet{wets1985continuity}, \citet{zhang1990note}), quadratic (\citet{lee2006continuity}), and quasiconvex (\citet{terazono2015continuity}). 
Also, there are generalized versions of Maximum Theorem (\citet{walker1979generalization}, \citet{leininger1984generalization}, \citet{ausubel1993generalized}). 
We handle these analyses to bridge parametric optimization and NN.

Also, parallel to the development of various optimization methodologies for ML (\citet{wright2022optimization}, \citet{boyd2004convex}, \citet{sra2012optimization}), there has been increasing interest both from operations research and computer science communities to solve mathematical optimization problems using ML. 
The literature on learning parametric optimization shows two main approaches. 
The first approach is to learn the optimal solution directly from the input parameter by utilizing the existing solution as data (\citet{lillicrap2015continuous}, \citet{mnih2015human}, \citet{vinyals2015pointer}, \citet{dai2017learning}, \citet{li2016learning}, \citet{donti2017task}). 
Although this approach is simple and intuitive, it has the disadvantage of providing solutions that may violate the critical constraints of optimization problems. 
This prompted second approach that applies ML indirectly as intermediate step. \citet{li2018combinatorial} used Graphical Neural Networks to guide a parallelized tree search procedure that rapidly generate a large number of candidate solutions. \citet{agrawal2020learning} applied ML to approximate gradient of solution of convex optimization. \citet{misra2021learning} and \citet{bertsimas2019online} identified optimal active or tight constraint sets. 
For more recent works, see \citet{bae2023deep, dai2021neural, dumouchelle2022neur2sp, pmlr-v202-kim23r}.
In this paper, we consider the situation where the parametric optimization problems are approximated by NN directly for utilization of UAT.
\subsection{Contribution.}\label{cont}
In this paper, we derive conditions for UAT to hold for approximating parametric optimization problems. With our derivations, we can specify how to formulate the parametric optimization problem rather than naively hoping that NN will approximate the optimization problem well. The main contribution of the study can be summarized as below.
\begin{itemize}
    \item We provide sufficient conditions for UAT to hold for optimal policy for continuous parametric optimization problems and these conditions are quite general in that we do not impose convexity or even quasi-convexity of optimization problems. We only exploit the assumptions and results of the Maximum Theorem (\citet{berge1963espaces}).
    \item We also address situations when these sufficient conditions are not satisfied. In particular, we define a sampling function and its stability which makes good approximation possible even without the sufficient conditions in original problems. Under the stable sampling function, original problems become reduced problem in which all conditions in main theorem are satisfied.
    \item We directly link vast amount of literature on NN with approximating optimization problems. There are many literatures linking the specific structure of NN to UAT. However, to our best knowledge, the general connection between the structure of parametric optimization problem and UAT has been scarcely investigated from the theoretical point of view. Our research clarifies such a vague connection by constructing piecewise linear policy approximation for NN.    
\end{itemize}

\subsection{Outline.}\label{outline}
The remainder of the paper is organized as follows. Preliminaries for deriving our results are included in Section \ref{preli}, and our main results with piecewise linear policy approximation are presented in Section \ref{main}. Section \ref{approxnn} discuss suitability of NN as estimator of our policy approximation. Improving feasibility and training with suboptimal training data are discussed in Section \ref{practical}. Finally, Section \ref{conclu} concludes the paper.

\section{Preliminaries.}\label{preli}
In Section \ref{preli}, we begin by introducing definitions and notations that are necessary for deriving our results. 
We also formally define the problem and list its assumptions.
\subsection{Definitions and Notations.}\label{defnot}
Parametric optimization takes the form,
\[
\min_x\quad 
f(x,\theta) \tag{1} \quad \text{subject\;to} \quad x\in C(\theta). \label{eq:1} 
\]
where $x\in X\subset \mathbb{R}^n$ is the decision variable, $\theta\in\Theta\subset\mathbb{R}^k$ is the parameter, $f:\mathbb{R}^n\times\mathbb{R}^k\rightarrow\mathbb{R}$ is the objective function and $C: \mathbb{R}^k \rightrightarrows 2^{\mathbb{R}^n}$ is a multivalued mapping, or correspondence, representing the feasible region defined by a set of constraints parameterized by $\theta$. 
Let the optimal value function $f^\ast:\mathbb{R}^k\rightarrow\mathbb{R}$ by $f^\ast(\theta)=\min_{x \in C(\theta)} f(x,\theta)$. 
We denote the optimal policy correspondence $C^\ast:\mathbb{R}^k\rightrightarrows2^{\mathbb{R}^n}$ by $C^\ast(\theta)=\arg\min_{x \in C(\theta)}f(x,\theta)=\{x\in C(\theta)|f(x,\theta)=f^\ast(\theta)\}$. 
An optimal solution $x^\ast(\theta)$ is an element of $C^\ast(\theta)$.

For any vector $x\in\mathbb{R}^n$, its norm $\| x\|$ is defined as the Euclidean norm, $\| x\|^2=\sum_{i=1}^{n} x^2_i$. 
For any non-empty set $X$ of vectors in $\mathbb{R}^n$, the $\varepsilon$-neighborhood is represented by $\mathcal{B}_\varepsilon(X)=\{y\in \mathbb{R}^n|$  $\exists x\in X$ s.t. $\| x-y\|<\varepsilon\}$. 
We define the stability of correspondence based on the continuity by \citet{hausdorff2021set}.
While there are different definitions of stability (\citet{berge1963espaces}, \citet{hogan1973point}), the Hausdorff’s version is the most general (\citet{zhao1997lower}).
\begin{definition}
Let $C$ be a correspondence from parameter space $\mathrm{\Theta}\subset\mathbb{R}^k$ to $2^{\mathbb{R}^n}$. Then,
\begin{enumerate}
    \item[(a)] $C$ is upper hemicontinuous at $\theta_0$ if $\forall \varepsilon>0$,  $\exists \delta>0$ s.t. $C(\theta)\subset\mathcal{B}_\varepsilon(C(\theta_0)), \forall\theta\in\mathcal{B}_\delta(\theta_0)$
    \item[(b)] $C$ is lower hemicontinuous at $\theta_0$ if $\forall \varepsilon>0$,  $\exists \delta>0$ s.t. $C(\theta_0)\subset\mathcal{B}_\varepsilon(C(\theta)), \forall\theta\in\mathcal{B}_\delta(\theta_0)$
    \item[(c)] $C$ is continuous at $\theta_0$ if $C$ is both upper and lower hemicontinuous at $\theta_0$    
\end{enumerate}
\label{def_cor}
\end{definition}

UAT describes the capability of NN as an approximator. 
Although there are many variations, the key statement is that a function expressed as NN is dense on the function space of interest. 
The most classical version of UAT is independently introduced by \citet{hornik1989multilayer}. 
Since we are utilizing the key findings of UAT, we summarize and restate this study as presented in Theorem \ref{thm_uni}, where the term \textit{function} is written as a \textit{single-valued function} to distinguish it from a correspondence.

\begin{theorem}[Universal Approximation Theorem, restated from \citet{hornik1989multilayer}]
Let $f$ be a continuous single-valued function on a compact set $K$. Then, there exists a feed forward NN with a single hidden layer that uniformly approximates $f$ to within an arbitrarily $\varepsilon>0$ on $K$.
\label{thm_uni}
\end{theorem}

\subsection{The Maximum Theorem.}
The Maximum Theorem was presented in \citet{berge1963espaces}, which provides conditions under which the value function is continuous and the optimal policy correspondence is upper hemicontinuous for a parametric optimization problem given by (\ref{eq:1}). 
This theorem sets the basis for developing a connection between parametric optimization and UAT. 
We restate Berge’s Maximum Theorem as Theorem \ref{thm_max}.
\begin{theorem}[The Maximum Theorem, restated from \citet{berge1963espaces}]
Let $f\mathrm{\,:\,}X\times\mathrm{\Theta}\rightarrow\mathbb{R}$ be a continuous function on the product $X\times\mathrm{\Theta}$, and $C\mathrm{\,:\,} \Theta\rightrightarrows X$ be a compact-valued correspondence s.t. $C(\theta)\neq\emptyset, \forall\theta\in\mathrm{\Theta}$. Define the $f^\ast(\theta)$ and $C^\ast(\theta)$ as $f^\ast(\theta)=\min_{x \in C(\theta)} f(x,\theta)$ and $C^\ast(\theta)=\arg\min_{x \in C(\theta)} f(x,\theta)=\{x\in C(\theta)|f(x,\theta)=f^\ast(\theta)\}$. If $C$ is continuous (i.e. both upper and lower hemicontinuous) at $\theta$, then $f^\ast$ is continuous and $C^\ast$ is upper hemicontinuous with non-empty and compact values.
\label{thm_max}
\end{theorem}

\subsection{Problem Description.}\label{probdef}
Our goal is to find the conditions of $f$ and $C$ that allows UAT to be applied to approximating the optimal policy correspondence $C^\ast$. 
Suppose the optimization problem given by (\ref{eq:1}) is formulated so that it changes stably as $\theta$ varies. The key questions are as follows.
\begin{enumerate}
\item[(Q1)] Is $C^\ast$ continuous or single-valued function?
\item[(Q2)] Are there bounds on errors from approximation, and do they converge to zero?
\item[(Q3)] Is NN suitable class for learning of $C^\ast$? 
\end{enumerate}

Questions (Q1) arise as UAT generally requires continuity and a single-valued function. 
We analyze (Q1) based on the Maximum Theorem (\citet{berge1963espaces}), which is one of the most applied theorems in stability theory. 
To guarantee an acceptable approximation, we construct a target function for optimal policy $C^\ast$, which is a piecewise linear continuous function and derive conditions where the approximation error converges to zero. 
This will address question (Q2). 
Finally, for question (Q3), we represent generalization error and approximation error of NN on learning constructed piecewise linear continuous target function.

\subsection{Assumptions.}\label{assum}
For problem (1), we assume that the objective function $f(x,\theta)$ is continuous on the product $X\times\Theta$, the feasible region $C(\theta)$ is continuous on $\Theta$ and $C(\theta)$ is a non-empty compact set for each $\theta\in\Theta$.
We make assumptions on the training data for optimal policy as well. 
A training example for parametric optimization is a pair of a parameter and its corresponding optimal solution $(\theta_i,x^\ast(\theta_i))$. 
Let the training data be the set of examples, $T=\{(\theta_i,x^\ast(\theta_i))|i=1,\cdots,m\}$. 
Notice that there can be more than one optimal solution $x^\ast(\theta_i)$ for each $\theta_i$. 
In practice, it is computationally expensive, if not impossible, to obtain the entire set of optimal solutions.
In fact, it is difficult even to identify whether there are multiple optimal solutions or not. 
Therefore, to incorporate such practical aspects, we assume that there exists a solver that can extract exactly one element $x^\ast(\theta_i)$ from $C^\ast(\theta_i)$ for any given $\theta_i$. 
However, it does not have control on the choice of $x^\ast(\theta_i)$, so that the optimal solution is obtained in a random manner from $C^\ast(\theta_i)$. 
Moreover, the solver is not able to identify if $C^\ast(\theta_i)$ is a singleton or not. 
It is as if the training data is a discrete sample path from the correspondence $C^\ast$ indexed by $\{\theta_i|i=1,\cdots,m\}$.

\section{Piecewise Linear Policy Approximation.}\label{main}
In Section \ref{main}, we present our main results about piecewise linear policy approximation.
Given the above assumptions on $f$ and $C$, the Theorem \ref{thm_max} states that $f^\ast$ is a continuous function, and $C^\ast$ is a non-empty and compact-valued upper hemicontinuous correspondence. Thus, unlike the value function $f^\ast$, which guarantees universal approximation, $C^\ast$ is not a single-valued function and is not even continuous, which requires additional treatments. Before making further steps, we state the following as a special case.

\begin{corollary}
If $C^\ast$ is singleton for each $\theta$, NN universally approximates $C^\ast$.
\label{cor_sing}
\end{corollary}

If the optimal solution for (\ref{eq:1}) is unique for every $\theta$, its optimal policy is not a correspondence, and reduces to a single-valued function. 
As upper hemicontinuity implies the continuity for a function, UAT can readily be applied for $C^\ast$. 
While this is a special case with a strong assumption, Corollary \ref{cor_sing} is the ideal case. 
In general, there can be multiple optimal solutions for some $\theta$, and, thus, $C^\ast$ is no longer a single-valued function. 
But under some conditions on $C^\ast$, there is possibility to find a continuous function called a \textit{selection} as defined in Definition \ref{def_sel}. 

\begin{definition}
Given two sets $X$ and $Y$, let $F$ be a correspondence from $X$ to $Y$. A function $f\mathrm{\,:\,}X\rightarrow Y$ is said to be a selection of $F$, if \;$\forall x\in X, f(x)\in F(x)$.
\label{def_sel}
\end{definition}

\begin{proposition}[Existence of a continuous selection]
$C^\ast$ has a continuous selection if it is convex-valued and lower hemicontinuous.
\label{pro_sel}
\end{proposition}
\proof
This is an immediate result of the selection theorem by Michael \cite{michael1956continuous}.

\begin{figure}[t]
\centering
\includegraphics[scale=0.4]{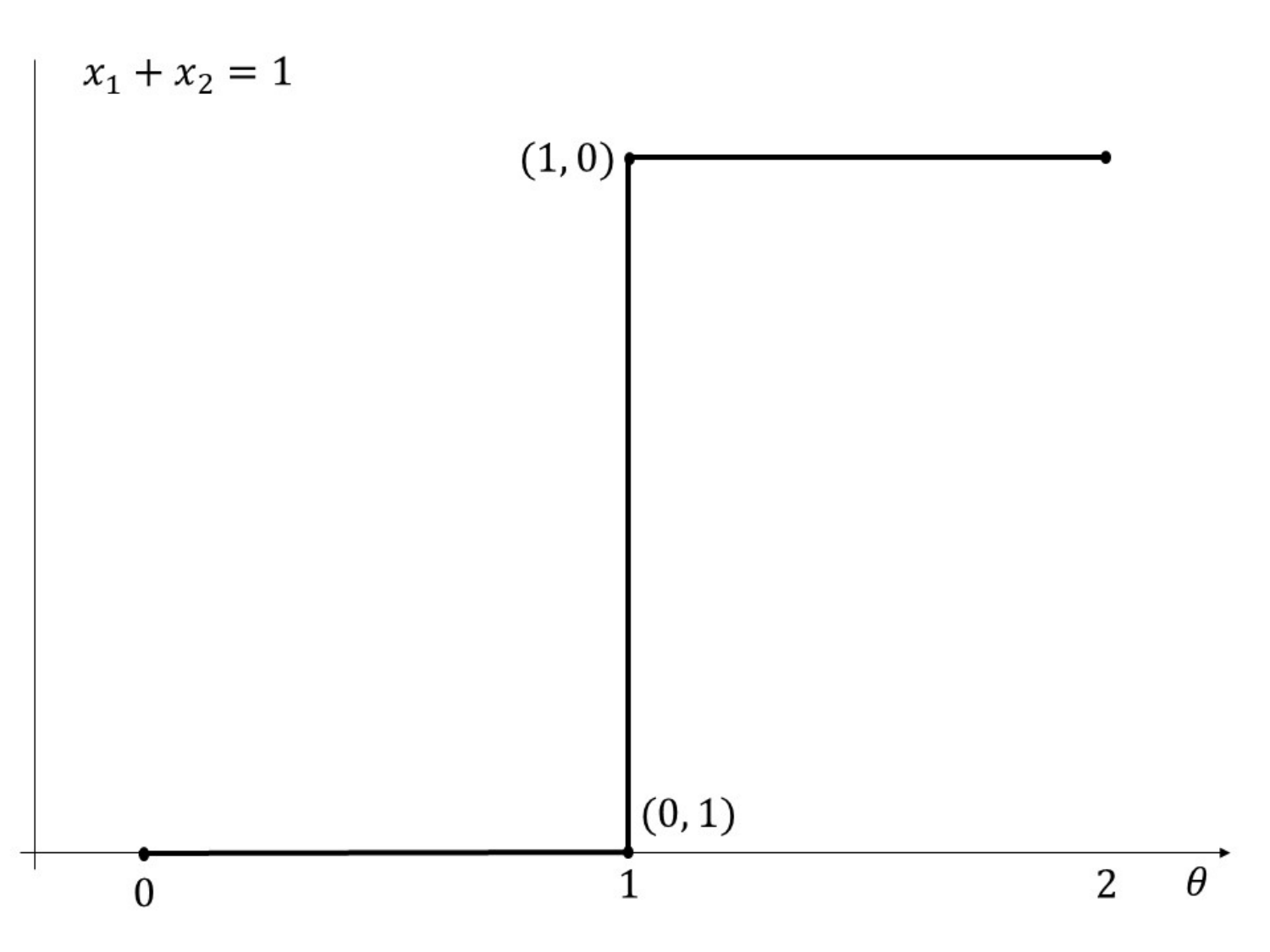}
\caption{Optimal policy correspondence of Example \ref{ex1} on $x_1+x_2=1$}
\label{fig:1}
\end{figure}

\endproof
There is a potential issue with Proposition \ref{pro_sel}. 
Some important classes of optimization, including linear programming problems, do not necessarily have lower hemicontinuous optimal policy correspondence. 
To illustrate the issues on approximating $C^\ast$, consider the following linear program with a parameter $\theta\in[0,2]$,
\begin{example}
\[
\min_{x_1,x_2} \quad -\theta x_1 - x_2 \quad \text{subject\;to} \quad x_1 + x_2 \leq 1, \quad  x_1 \geq 0, \quad x_2 \geq 0.
\]
\label{ex1}
\end{example}
The optimal policy correspondence for the problem given by Example \ref{ex1} becomes
\[
C^\ast(\theta) = \begin{cases}
\{(0,1)\} & \text{for $\theta \in [0, 1)$},\\
\{(x_1,x_2)| x_1+x_2=1,\ x_1\geq0,x_2\geq0\} & \text{for $\theta=1$},\\
\{(1,0)\}\ & \text{for $\theta \in (1, 2]$.}
\end{cases}
\]
As illustrated in Figure \ref{fig:1}, $C^\ast(\theta)$ contains a jump at $\theta=1$. 
Thus, it is evident that there is no continuous selection of $C^\ast(\theta)$. 
This means that UAT cannot be directly applied, and we need to find a workaround to make it work. 

Thus, we drop the assumption that $C^\ast$ is lower hemicontinuous and only take upper hemicontinuity of $C^\ast$, which is guaranteed by Theorem \ref{thm_max}. 
Since a continuous selection generally does not exist, we construct a new target function.
For given training data $T=\{(\theta_i,x^\ast(\theta_i))|i=1,\cdots,m\}$, we attempt to estimate the optimal solution on the convex hull of $\theta_1,\theta_2,\ldots,\theta_m$ denoted as $Conv(\{\theta_1,\theta_2,\ldots,\theta_m\})$. 
Furthermore, we consider a finite collection $\mathcal{S}=\{S_1,\ldots,S_d\}$ of subset of $\{\theta_1,\theta_2,\ldots,\theta_m\}\subset\mathbb{R}^k$ such that
\begin{enumerate}
    \item[(a)] For each $j\in\{1,\ldots,d\}$, $|S_j|=k+1$ i.e. $Conv(S_j)$ is $k$-dimensional simplex.
    \item[(b)] $Conv(\{\theta_1,\theta_2,\ldots,\theta_m\})=\bigcup_{j=1}^dConv(S_j)$.
    \item[(c)] For any non-empty subset $J \subseteq \{1,2,\ldots,d\}$, $\bigcap_{j \in J}Conv(S_j)=Conv(\bigcap_{j \in J}S_j)$.
\end{enumerate}

This collection is called triangulations (\citet{lee2017subdivisions}).
One way of constructing this is lexicographic triangulations introduced by \citet{sturmfels1991grobner}. Given such collection $\mathcal{S}$, for any $\theta \in Conv(\{\theta_1,\theta_2,\ldots,\theta_m\})$, there exists index $j_\theta \in \{1,\ldots,d\}$ such that  $\theta=\lambda^{j_\theta}_1\theta^{j_\theta}_1+\ldots+\lambda^{j_\theta}_{k+1}\theta^{j_\theta}_{k+1}$ where $\{\theta^{j_\theta}_1,\ldots,\theta^{j_\theta}_{k+1}\}=S_{j_\theta}\in \mathcal{S}$ and $\lambda^{j_\theta}_1,\ldots,\lambda^{j_\theta}_{k+1}\geq0,\lambda^{j_\theta}_1+\ldots+\lambda^{j_\theta}_{k+1}=1$. 
Our key approach is to approximate $x^\ast(\theta)$ as $\lambda^{j_\theta}_1{x^\ast(\theta}_1^{j_\theta})+\ldots+\lambda^{j_\theta}_{k+1}{x^\ast(\theta}_{k+1}^{j_\theta})$. This approximation is a single-valued by construction, and continuous function 
by following theorem, so that UAT can be applied.

\begin{theorem}[Continuity of Policy Approximation]
Consider a finite collection $\mathcal{S}=\{S_1,\ldots,S_d\}$ of subset of $\{\theta_1,\theta_2,\ldots,\theta_m\}\subset\mathbb{R}^k$ such that
\begin{enumerate}
    \item[(a)] For each $j\in\{1,\ldots,d\}$, $|S_j|=k+1$ i.e. $Conv(S_j)$ is $k$-dimensional simplex.
    \item[(b)] $Conv(\{\theta_1,\theta_2,\ldots,\theta_m\})=\bigcup_{j=1}^dConv(S_j)$.
    \item[(c)] For any non-empty subset $J \subseteq \{1,2,\ldots,d\}$, $\bigcap_{j \in J}Conv(S_j)=Conv(\bigcap_{j \in J}S_j)$.
\end{enumerate}  
Denote $S_j=\{\theta^j_1,\ldots,\theta^j_{k+1}\}$. For any $\theta\in Conv(\{\theta_1,\theta_2,\ldots,\theta_m\})$,  let $j_\theta$ be index of $\mathcal{S}$ such that $\theta \in Conv(S_{j_\theta})$. Then, the function $\hat{x}:Conv(\{\theta_1,\theta_2,\ldots,\theta_m\}) \rightarrow \mathbb{R}^n$
\[
\hat{x}(\theta) = \lambda^{j_\theta}_1{x^\ast(\theta}_1^{j_\theta})+\ldots+\lambda^{j_\theta}_{k+1}{x^\ast(\theta}_{k+1}^{j_\theta}).
\]
is continuous where $\lambda^{j_\theta}_1,\ldots,\lambda^{j_\theta}_{k+1}$ are weights of convex combination for $\theta$ i.e. $\theta=\lambda^{j_\theta}_1\theta^{j_\theta}_1+\ldots+\lambda^{j_\theta}_{k+1}\theta^{j_\theta}_{k+1}$
\label{thm_con}
\end{theorem}
\proof
We first prove Lemma \ref{lem_con}.
\begin{lemma}
Let $X_1, X_2, \ldots, X_d$ be closed subsets of  $\mathbb{R}^k$. Suppose that  $f_j : X_j \rightarrow \mathbb{R}^n, j=1, \ldots, d$  are continuous functions such that for any nonempty subset $J \subseteq \{1,2,\ldots,d\}$ and for any $j_1, j_2 \in J$ 
\[
f_{j_1}|\bigcap_{j\in J}X_j = f_{j_2}|\bigcap_{j\in J}X_j 
\]
Then, the function
\[
f : \bigcup_{j=1}^d X_j \rightarrow \mathbb{R}^n, \quad x \rightarrow
\begin{cases}
       f_1(x), \quad x \in X_1 \\
       f_2(x), \quad x \in X_2 \\
       \vdots \\
       f_d(x), \quad x \in X_d \\
\end{cases}
\]
is also continuous. 
\label{lem_con}
\end{lemma}
\proof[Proof of Lemma \ref{lem_con}]
Let $A \subseteq \mathbb{R}^n$ be closed in $\mathbb{R}^n$. Because $f_j$ is continuous, $f_j^{-1}(A)$ is closed in $X_j$. More formally, there is $Y_j \subseteq \mathbb{R}^k$ such that $Y_j$ is closed in $\mathbb{R}^k$ and $f^{-1}(A)=Y_j \cap X_j$.  Then, we have
\[
f^{-1}(A)=\bigcup_{j=1}^d f_j^{-1}(A)=\bigcup_{j=1}^d(Y_j \cap X_j)=(X_1\cup\ldots\cup X_d)\cap(\bigcap_{Z_j\in\{X_j, Y_j\}, j=2,\ldots,d} (Y_1\cup Z_2 \cup\ldots\cup Z_d)).
\]
If follows from $Y_j, X_j$ are closed in $\mathbb{R}^k$ that $B:=\bigcap_{Z_j\in\{X_j, Y_j\}, j=2,\ldots,d} (Y_1\cup Z_2 \cup\ldots\cup Z_d)$ are closed in $\mathbb{R}^k$. Hence, $f^{-1}(A)=(X_1\cup\ldots\cup X_d)\cap B$ is closed in  $(X_1\cup\ldots\cup X_d)$. Thus, $f$ is continuous.
\endproof
Now we prove Theorem \ref{thm_con}. Define function $\hat{x}_j:Conv(S_j) \rightarrow \mathbb{R}^n$ as $\hat{x}_j(\theta)=\lambda^j_1{x^\ast(\theta}_1^j)+\ldots+\lambda^j_{k+1}{x^\ast(\theta}_{k+1}^j)$ for $\theta = \lambda_1^j \theta_1^j+\ldots+\lambda_{k+1}^j \theta_{k+1}^j$. Next, we prove the function $\hat{x}_j$ is continuous.  For $\theta = \lambda_1^j \theta_1^j+\ldots+\lambda_{k+1}^j \theta_{k+1}^j \in Conv(S_j)$, let $y_j(\theta) = (\lambda^j_1, \ldots, \lambda^j_{k+1})$. Then, $y_j(\theta)$ is inverse of linear function (which is linear) and $\hat{x}_j(\theta)=y_j(\theta)^T(x^\ast(\theta^j_1),\ldots,x^\ast(\theta^j_{k+1})) $ is linear function of $y_j(\theta)$. Thus, $\hat{x}_j(\theta)$ is composite function of two linear functions which is continuous. 
Also, for any nonempty subset $J \subseteq \{1,2,\ldots,d\}$ and for any $j_1, j_2 \in J$,
\[
\hat{x}_{j_1}|\bigcap_{j\in J}Conv(S_j) = \hat{x}_{j_2}|\bigcap_{j\in J}Conv(S_j).
\]
since $\bigcap_{j \in J}Conv(S_j)=Conv(\bigcap_{j \in J}S_j)$. Thus, by the Lemma \ref{lem_con}, the function $\hat{x}$
\[
   \hat{x} : Conv(\{\theta_1,\theta_2,\ldots,\theta_m\}) \rightarrow \mathbb{R}^n, \quad \theta \rightarrow
   \begin{cases}
       \hat{x}_1(\theta), \quad \theta \in Conv(S_1) \\
       \hat{x}_2(\theta), \quad \theta \in Conv(S_2) \\
       \vdots \\
       \hat{x}_d(\theta), \quad \theta \in Conv(S_d).
   \end{cases}
\]
is continuous function.
\endproof

An inherent question regarding the function $\hat{x}$ pertains to the degree of approximation accuracy.  We first remark that there exists parametric optimization problem where convergence of errors of $\hat{x}$ is not guaranteed since $x^\ast(\theta)$ is arbitrarily chosen from $C^\ast(\theta)$. 
\begin{example}
\[
\min_x\quad f(x,\theta)=(x-1)^2(x-3)^2 \quad \text{subject\;to} \quad x\in C(\theta)=\{1,3\}.
\]
\label{ex2}
\end{example}
For the trivial parametric optimization problem given by Example \ref{ex2},
the optimal policy correspondence is
\[
C^\ast(\theta)=\{1,3\}, \forall\theta
\]
In this case, our construction always has suboptimality of 1 for some $\theta$ because $|f^\ast({(\theta}_1+\theta_2)/2)-f(2,\ {(\theta}_1+\theta_2)/2)|=1$ if $(\theta_1, 1)$ and $(\theta_2, 3)$ are sampled with $\theta_1<\theta_2$. 

In order to determine the suitability of using $\hat{x}$ for approximating specific parametric optimization problems, we establish metrics to assess the performance of the target function $\hat{x}$. These metrics are referred to as $\varepsilon$-\textit{suboptimality} and $\varepsilon$-\textit{infeasibility}.
We present our main development in Theorem \ref{thm_main}, which states that a constructed function is $\varepsilon$-\textit{infeasible} and $\varepsilon$-\textit{suboptimal} solution for sufficiently dense training data under certain conditions. 
While this adds some restrictions, it allows applying UAT to parametric optimization and these two conditions can be lifted with a stable sampler, which we further discuss in the subsequent sections.

\begin{definition}
Let $f$ be the objective function and $C$ be the correspondence in the formulation given by (\ref{eq:1}). Then, for $\varepsilon>0$,
\begin{enumerate}
    \item[(a)] $\hat{x}(\theta)$ is $\varepsilon$-infeasible solution if $\hat{x}(\theta)\in\mathcal{B}_{\varepsilon}(C(\theta)).$
    \item[(b)] $\hat{x}(\theta)$ is $\varepsilon$-suboptimal solution if $|f(\hat{x}(\theta),\theta)-f(x^\ast(\theta),\theta)|<\varepsilon.$
\end{enumerate}
\label{def_eps}
\end{definition}

\begin{theorem}[Convergence Property of Piecewise Linear Policy Approximation]
Suppose that $f$, $C$ satisfy all conditions for Theorem \ref{thm_max}. 
Define a training data set $T=\{(\theta_i,x^\ast(\theta_i))|\theta_i\in\mathrm{\Theta},i=1,\cdots,m\}$ where $x^\ast(\theta_i)$ is an arbitrarily chosen point from $C^\ast(\theta_i)$.
For $\theta\in Conv(\{\theta_1,\theta_2,\ldots,\theta_m\})$, define $d(S_{j_\theta})=\max\{\| p-q\|:p,q\in Conv(S_{j_\theta})\}$ where $S_{j_\theta}$ is element of finite collection $\mathcal{S}$ in Theorem \ref{thm_con}. Then, the function $\hat{x}(\theta)$ in Theorem \ref{thm_con} satisfies the followings.
\begin{enumerate}
    \item[(a)] If $Conv(C^\ast(\theta))\subseteq C(\theta)$, for given $\varepsilon>0$, $\exists\delta>0$ s.t. if $d(S_j)<\delta$, $\hat{x}(\theta)$ is $\varepsilon$-infeasible solution.
    \item[(b)] If $f(x,\theta)=f^\ast(\theta), \forall x\in Conv(C^\ast(\theta))$, for given $\varepsilon>0$, $\exists\delta>0$ s.t. if $d(S_j)<\delta$, $\hat{x}(\theta)$ is $\varepsilon$-suboptimal solution.
\end{enumerate}
\label{thm_main}
\end{theorem}
\proof
(a) Assume $Conv(C^\ast(\theta))\subseteq C(\theta)$ and let $\varepsilon>0$. 
Since $f$, $C$ satisfy all conditions for Theorem \ref{thm_max}, $C^\ast$ is upper hemicontinuous. 
Thus, a set 
\[
\mathrm{\Delta}_\theta=\{\delta>0:\|\theta-\theta^\prime\|< \delta\Longrightarrow C^\ast(\theta^\prime)\subseteq\mathcal{B}_\varepsilon(C^\ast(\theta))\}.
\]
is not empty.

Define $\delta_\theta=\sup{\mathrm{\Delta}_\theta}$. 
Choose $\delta=\delta_\theta$. 
If $d(S_{j_\theta})<\delta$, $\|\theta-\theta^{j_\theta}_l\|<\delta$ i.e. $C^\ast(\theta_l^{j_\theta})\subseteq\mathcal{B}_\varepsilon(C^\ast(\theta)), \forall l=1,\ldots,k+1$.
Then, with the assumption and $C^\ast(\theta)\subseteq Conv(C^\ast(\theta))$,
\[
x^\ast(\theta_l^{j_\theta})\in C^\ast(\theta_l^{j_\theta})\subseteq\mathcal{B}_\varepsilon(C^\ast(\theta))\subseteq\mathcal{B}_\varepsilon(Conv(C^\ast(\theta)))\subseteq\mathcal{B}_\varepsilon(C(\theta)),\forall l=1,2,\ldots,k+1.
\]
Thus, $x^\ast(\theta_l^{j_\theta})\in\mathcal{B}_\varepsilon(Conv(C^\ast(\theta)))\subseteq\mathcal{B}_\varepsilon(C(\theta))$. 
Note that, since $Conv(C^\ast(\theta))$ is convex set, $\mathcal{B}_\varepsilon(Conv(C^\ast(\theta)))$ is also convex set. 
This means the convex combination 
\[
\lambda^{j_\theta}_1{x^\ast(\theta}_1^{j_\theta})+\ldots+\lambda^{j_\theta}_{k+1}{x^\ast(\theta}_{k+1}^{j_\theta}).
\]
is in $\mathcal{B}_\varepsilon(Conv(C^\ast(\theta)))$. 
Thus, $\hat{x}(\theta)\in B_\varepsilon(C(\theta))$.

(b) We first prove Lemma \ref{lem_main}.
\begin{lemma}
If $A\subseteq\mathbb{R}^n$ is a compact set, $Conv(A)$ is also a compact set.
\label{lem_main}
\end{lemma}
\proof[Proof of Lemma \ref{lem_main}]
Carathéodory Theorem (Danninger-Uchida \cite{Danninger-Uchida2009}) states that each element of the convex hull of $A$ is a convex combination of $n+1$ elements of $A$. 
By defining $Simp(n)=\{(w_0,\ldots,w_n): w_j\geq0,w_0+\ldots+w_n=1\}$ and $F(a_0,\ldots,a_n;\ w_0,\ldots,w_n)=\sum_{i} w_ia_i$, $Conv(A)$ can be expressed as the image of the compact set $A^{n+1}\times Simp(n)$ under a continuous map $F$, and so it is compact.
\endproof
Now assume $f(x,\theta)=f^\ast(\theta), \forall x\in Conv(C^\ast(\theta))$ and let $\varepsilon>0$. 
We first show there exists $\delta^\prime>0$ such that, 
\[
\inf_{x^\ast\in Con v(C^\ast(\theta))} \|x^\prime-x^\ast\|< \delta^\prime\Longrightarrow |f(x^\prime,\theta)-f^\ast(\theta)|<\varepsilon, \forall x^\prime\in\mathcal{B}_\varepsilon(C(\theta)).
\]
Since $C$ is compact-valued correspondence, $\mathcal{B}_\varepsilon(C(\theta))\times{\theta}$ is compact set.
Thus, $f$ is uniformly continuous on $\mathcal{B}_\varepsilon(C(\theta))\times{\theta}$. 
Thus, there exist $\delta^{\prime\prime}>0$ such that, 
\[
\|y-z\|< \delta^{\prime\prime}\Longrightarrow |f(y,\theta)-f(z,\theta)|<\varepsilon, \forall y,z\in\mathcal{B}_\varepsilon(C(\theta)).
\]
Choose $\delta^\prime=\delta^{\prime\prime}$. 
Note that $C^\ast$ is compact-valued correspondence from Theorem \ref{thm_max}. 
Thus, $Conv(C^\ast(\theta))$ is also compact set from Lemma \ref{lem_main}. 
Hence, 
\[
x_{min}^\ast={\arg\min}_{x^\ast\in Conv(C^\ast(\theta))}\|x^\prime-x^\ast\|={\arg\min}_{x^\ast\in Conv(C^\ast(\theta))}\|x^\prime-x^\ast\|.
\] 
is in $Conv(C^\ast(\theta))$.
Now we have
\[
\inf_{x^\ast\in Conv(C^\ast(\theta))} \|x^\prime-x^\ast\|<\delta^\prime\Longrightarrow \|x^\prime-x_{min}^\ast\|<\delta^\prime=\delta^{\prime\prime}\Longrightarrow |f(x^\prime,\theta)-f(x_{min}^\ast,\theta)|<\varepsilon.
\]
Since $x_{min}^\ast\in Conv(C^\ast(\theta))$, $f(x_{min}^\ast,\theta)=f^\ast(\theta)$. 
Thus, $|f(x^\prime,\theta)-f^\ast(\theta)|<\varepsilon$

Now, we prove part (b) of the theorem. 
From the above statement, there exists $\delta^\prime>0$ such that, 
\[\inf_{x^\ast\in Con v(C^\ast(\theta))} \|x^\prime-x^\ast\|<\delta^\prime\Longrightarrow |f(x^\prime,\theta)-f^\ast(\theta)|<\varepsilon, \forall x^\prime\in\mathcal{B}_\varepsilon(C(\theta)).\]
Also, since $f$, $C$ satisfy all conditions for Theorem \ref{thm_max}, $C^\ast$ is upper hemicontinuous. 
Thus, a set
\[
\mathrm{\Delta}_\theta=\{\delta>0|\;if\;\|\theta-\theta^\prime\|< \delta,C^\ast(\theta^\prime)\subseteq\mathcal{B}_\delta^\prime(C^\ast(\theta))\}.
\]
is not empty.

Define $\delta_\theta=\sup{\mathrm{\Delta}}_\theta$.
Choose $\delta=\delta_\theta$. 
If $d(S_j)<\delta$, $\|\theta-\theta^{j_\theta}_l\|< \delta$ i.e. $C^\ast(\theta_l^{j_\theta})\subseteq\mathcal{B}_{\delta^\prime}(C^\ast(\theta)),\forall l=1,\ldots,k+1$.
Then with $C^\ast(\theta)\subseteq Conv(C^\ast(\theta))$,
\[
x^\ast(\theta_l^{j_\theta})\in C^\ast(\theta_l^{j_\theta})\subseteq\mathcal{B}_{\delta^\prime}(C^\ast(\theta))\subseteq\mathcal{B}_{\delta^\prime}(Conv(C^\ast(\theta))), \forall l=1,\ldots,k+1.
\]
Since $Conv(C^\ast(\theta))$ is convex set, $\mathcal{B}_{\delta^\prime}(Conv(C^\ast(\theta)))$ is also convex set. 
This means the convex combination 
\[
\hat{x}(\theta)=\lambda^{j_\theta}_1{x^\ast(\theta}_1^{j_\theta})+\ldots+\lambda^{j_\theta}_{k+1}{x^\ast(\theta}_{k+1}^{j_\theta}).
\]
is in $\mathcal{B}_{\delta^\prime}(Conv(C^\ast(\theta)))$. 
Also, note that $\hat{x}(\theta)\in\mathcal{B}_\varepsilon(C(\theta))$ from part (a) since assumption in (b) indicates $Conv(C^\ast(\theta))\subseteq C^\ast(\theta)\subseteq C(\theta)$. 
Accordingly, $|f(\hat{x}(\theta),\theta)-f^\ast(\theta)|<\varepsilon$. 

\endproof
Theorem \ref{thm_main} shows that if the problem (\ref{eq:1}) satisfies the sufficient conditions, the errors on feasibility and optimality of our piecewise linear policy approximation $\hat{x}(\theta)$ converges to zero. For example, suppose that training data $T=\{(\theta_i,x^\ast(\theta_i))|\theta_i\in\mathrm{\Theta},\ i=1,\cdots,m\}$ is sampled from Example \ref{ex1}. With out loss of generality, suppose that $\theta_1 \leq \theta_2 \leq \ldots \leq \theta_m$. Then, the finite collection $\mathcal{S}$ in Theorem \ref{thm_con} can be constructed as $\mathcal{S}=\{S_j, j=1,\ldots,m-1\}$ where $S_j = \{\theta_j,\theta_{j+1} \}$.  Let $l$ be an index such that $\theta_l\le\ 1\le\theta_{l+1}$. Then  $\hat{x}(\theta)$ is
\[
\hat{x}(\theta) = \begin{cases}
(0,\ 1) & \text{$\theta\le\theta_l$},\\
(\frac{\theta-\theta_l}{\theta_{l+1}-\theta_l},\frac{\theta_{l+1}-\theta}{\theta_{l+1}-\theta_l}) & \text{$\theta_l<\theta<\theta_{l+1}$},\\
(1,\ 0) & \text{$\theta_{l+1}\le\theta$}.
\end{cases}
\]
Note that $\hat{x}(\theta)$ is a feasible solution, and, thus, an $\varepsilon$-infeasible solution. 
We want to further show that $\hat{x}(\theta)$ is an $\varepsilon$-suboptimal solution. 
It holds that $f^\ast(\theta)=f(\hat{x}(\theta),\theta)$ for the three cases: $\theta\le\theta_l, \theta\geq\theta_{l+1}$, or $\theta=1$. 
For $\theta_l<\theta<1$, $\varepsilon$-suboptimality can be shown as below if we choose $\delta<4\varepsilon$,
\[f^\ast(\theta)-f(\hat{x}(\theta),\theta)=1-\frac{\theta(\theta-\theta_l)}{\theta_{l+1}-\theta_l}-\frac{\theta_{l+1}-\theta}{\theta_{l+1}-\theta_l}\]
\[=\frac{(\theta-\theta_l)(1-\theta)}{\theta_{l+1}-\theta_l}<\frac{(\theta-\theta_l)(1-\theta)}{1-\theta_l}\le\frac{1-\theta_l}{4}<\frac{\delta}{4}<\varepsilon.\]
Similarly, the same results can be derived for $1<\theta<\theta_{l+1}$. 
Thus, $\hat{x}(\theta)$ is an $\varepsilon$-\textit{suboptimal solution}.

\subsection{General Case with a Stable Sampling Function.}\label{stable}
We have shown that the proposed piecewise linear policy approximation can be a reasonable solution under some assumptions on $f$ and $C$ for dense enough training data.
In this subsection, a stable sampling function designed to tackle a broader range of  parametric optimization problem.
Note that Example \ref{ex2} does not satisfy the conditions for Theorem \ref{thm_main} since $Conv(C^\ast(\theta))\nsubseteq C^\ast(\theta)$ and $Conv(C^\ast(\theta))\nsubseteq C(\theta)$. 
But even in this case, it may be possible to apply the Theorem \ref{thm_main} by sampling data from certain parts of $C^\ast(\theta)$, which has practical implications since decision makers often understand the nature of parametric optimization problems. 
Based on this idea, we define notion and stability of a sampling function.

\begin{definition}
Define a sampling function $s\mathrm{\,:\,}\mathrm{\Theta}\rightarrow\bigcup_{\theta}{C^\ast(\theta)}$ as $s(\theta)=x^\ast(\theta)$. 
Sampling function $s$ is stable with respect to $C^\ast$ if there exists a non-empty, compact, and convex-valued upper hemicontinuous correspondence $\overline{C^\ast}$ such that $s(\theta)\in\overline{C^\ast}(\theta)\subseteq C^\ast(\theta) \;\forall\theta.$
\label{def_sta}
\end{definition}

Note that the stable sampling function does not always exist. 
It depends on formulation of parametric optimization problem. 
For example, any sampling function in Example \ref{ex1} is stable since $C^\ast(\theta)$ is convex-valued. 
In Example \ref{ex2}, a sampling function that samples only 1 or 3 is stable, choose $\overline{C^\ast}$ as $\{1\}$ or $\{3\}$. 
Consider the following modified version of Example \ref{ex2}.
\begin{example}
\[
\min_x\quad f(x,\theta)=(x-1)^2(x-3)^2-\theta(x-3) \quad \text{subject\;to}\quad C(\theta)=\{1,3\}.
\]
\label{ex3}
\end{example}
The optimal policy correspondence for the problem given by Example \ref{ex3} is
\[
C^\ast(\theta) = \begin{cases}
\{1\} & \text{$\theta<0$},\\
\{1,\ 3\} & \text{$\theta=0$},\\
\{3\}\  & \text{$\theta>0$}
\end{cases}
\]
Note that there are two convex-valued sub correspondences of $C^\ast(\theta)$, $C^1(\theta)$ and $C^2(\theta)$, which neither is upper hemicontinuous.
\[
C^1(\theta) = \begin{cases}
\{1\} & \text{$\theta<0$},\\
\{1\} & \text{$\theta=0$},\\
\{3\}\  & \text{$\theta>0$}
\end{cases}
\quad C^2(\theta) = \begin{cases}
\{1\} & \text{$\theta<0$},\\
\{3\} & \text{$\theta=0$},\\
\{3\}\  & \text{$\theta>0$}
\end{cases}
\]

The advantage of a stable sampling function is that it makes Theorem \ref{thm_main} applicable even to parametric optimization problems that do not satisfy conditions $Conv(C^\ast(\theta))\subseteq C(\theta)$ and $f(x,\theta)=f^\ast(\theta), \forall x\in Conv(C^\ast(\theta))$. The following theorem shows that these two conditions become trivial if $C^\ast$ is convex-valued correspondence.
\begin{theorem}
If $C^\ast$ is convex-valued correspondence, the followings are hold.
\begin{enumerate}
    \item[(a)] $Conv(C^\ast(\theta))\subseteq C(\theta)$
    \item[(b)] $f(x,\theta)=f^\ast(\theta), \forall x\in Conv(C^\ast(\theta))$
\end{enumerate}
\label{thm_sta}
\end{theorem}
\proof
Since $C^\ast(\theta)$ is convex set, $Conv(C^\ast(\theta))=C^\ast(\theta)$. 
With this fact, we have

(a) $Conv(C^\ast(\theta))=C^\ast(\theta)\subseteq C^\ast(\theta)\subseteq C(\theta)$.

(b) $C^\ast(\theta)\subseteq C^\ast(\theta)$ implies $f(x,\theta)=f^\ast(\theta), \forall x\in Conv(C^\ast(\theta))$.

\endproof

If we have a stable sampling function $s$, $\overline{C^\ast}(\theta)$ can be considered as a substitute for $C^\ast(\theta)$. Since $\overline{C^\ast}(\theta)$ is convex-valued, the sufficient conditions $Conv(\overline{C^\ast}(\theta))\subseteq C(\theta)$ and $f(x,\theta)=f^\ast(\theta), \forall x\in Conv(\overline{C^\ast}(\theta))$ become redundant from Theorem \ref{thm_sta}. Furthermore, since $\overline{C^\ast}(\theta)$ is compact-valued and upper hemicontinuous, the same arguments in proof of Theorem \ref{thm_main} are applicable to $\overline{C^\ast}(\theta)$, so that we can apply Theorem \ref{thm_main} to more general parametric optimization.

\section{Suitability of Neural Networks.}\label{approxnn}
In previous section, we focused on constructing and evaluating the target function.
We build up our target function as continuous piecewise linear single-valued function and state its goodness in terms of $\varepsilon$-infeasibility and $\varepsilon$-suboptimality. In this section, we discuss capability of NN in approximating the target function. It is also possible to obtain a target function directly from training data. However, for a new parameter, this require to find $k+1$ points and its corresponding weights of the convex combination. This process requires considerable computation cost as the dimension of the parameter space increases, so it is not suitable for real-time decision making. On the other hand, the forward process of passing data to NN requires no additional processing time for finding the $k+1$ points once it is trained in advance. Then, is NN a really good estimator for our piecewise linear policy approximation? We answer this question by two aspect, generalization error (GE) and approximation error (AE).

In order to define and examine the GE and AE of NN, it is necessary to make choices regarding the architecture and loss function. 
This paper focuses on a specific type of neural network architecture, namely a feed forward network with rectified linear unit (ReLU) activations. 
The ReLU activation function is defined as  $ReLU(h) = \max(h, 0)$. 
The layers of a ReLU network, denoted as $h^l, \forall l=1,\ldots,L$, are defined recursively using the following relation: $h^l = ReLU(W^l h^{l-1} + b)$, where $h^0 = ReLU(b)$. Here, $W^l$ represents the weight matrices associated with each layer, and they satisfy the constraint that their infinity norms, denoted as $\|W^l\|_\infty$, are bounded by $B_l$. Additionally, we assume that $\|W^1\|$ is bounded by $B_1$. To sum up, the set of functions that can be represented by the output at depth $L$ in a ReLU network is denoted as $\mathcal{H}^L=\{h^L(\{W^l\}_{l=1}^L) : \|W^l\|_\infty \leq B_l, l \in [1,L], \|W^1\| \leq B_1\}$. 
For a loss function of NN, we prove following theorem
\begin{theorem}
Let $h(\theta)$ be ReLU NN and let $\hat{x}(\theta)$ be a constructed piecewise linear approximation.
\begin{enumerate}
    \item[(a)] $\forall \varepsilon>0, \exists \delta>0$ such that If $\|\hat{x}(\theta)-h(\theta)\|<\delta$ and $\hat{x}(\theta)$ is $\varepsilon/2$-infeasible solution, $h(\theta)$ is $\varepsilon$-infeasible solution. 
    \item[(b)] $\forall \varepsilon>0, \exists \delta>0$ such that If $\|\hat{x}(\theta)-h(\theta)\|<\delta$ and $\hat{x}(\theta)$ is $\varepsilon/2$-suboptimal solution, $h(\theta)$ is $\varepsilon$-suboptimal solution. 
\end{enumerate}
\label{thm_nn}
\end{theorem}
\proof
(a) Since $\hat{x}(\theta)$ is $\varepsilon/2$-infeasible solution, $\hat{x}(\theta)\in \mathcal{B}_{\varepsilon/2}(C(\theta))$. Choose $\delta=\varepsilon/2$. Then, since $\|\hat{x}(\theta)-h(\theta)\|<\delta$,
\[
h(\theta)\in\mathcal{B}_{\varepsilon/2}(\hat{x}(\theta))\subseteq\mathcal{B}_{\varepsilon}(C(\theta)).
\]
 Thus, $h(\theta)$ is $\varepsilon$-infeasible solution
\newline
(b) Since $\hat{x}(\theta)$ is $\varepsilon/2$-suboptimal solution, $|f(\hat{x}(\theta),\theta)-f(x^\ast(\theta),\theta)|<\varepsilon/2$
Also, since $f$ is continuous, there exists $\delta'>0$ such that, 
\[
\|\hat{x}(\theta)-h(\theta)\|<\delta'\Longrightarrow |f(\hat{x}(\theta),\theta)-f(h(\theta), \theta)|<\varepsilon/2
\]
Choose $\delta=\delta'$. Then, 
\[
|f(x^\ast(\theta),\theta)-f(h(\theta),\theta)|
\]
\[
=|f(x^\ast(\theta),\theta)-f(\hat{x}(\theta),\theta)+f(\hat{x}(\theta),\theta)-f(h(\theta),\theta)|
\]
\[
\leq |f(x^\ast(\theta),\theta)-f(\hat{x}(\theta),\theta)|+|f(\hat{x}(\theta),\theta)-f(h(\theta),\theta)|
\]
\[
<\varepsilon/2+\varepsilon/2=\varepsilon
\]
Thus, $h(\theta)$ is $\varepsilon$-suboptimal solution.

\endproof
Theorem \ref{thm_nn} says that if the NN is sufficiently close to our target function $\hat{x}(\theta)$, NN also gives $\varepsilon$-infeasible and $\varepsilon$-suboptimal solution. Thus, for a set of $m$ training examples  $T=\{(\theta_i, x^\ast(\theta_i))\}_{i=1}^m$ and ReLU NN $h$, we suppose that the weights are learned by minimizing a loss function $\ell$ which is defined as
\[
\ell(h(\theta), x^\ast(\theta))=\|h(\theta) - x^\ast(\theta)\|^2
\]
\subsection{Generalization error.}
For a hypothesis class $\mathcal{H}$, the GE is defined as 
\[
\mathcal{G}(\mathcal{H})=\mathbb{E}_T\sup_{h\in\mathcal{H}}(L_D(h)-L_T(h)).
\]
where $L_D(h)=\mathbb{E}_{\theta \sim D} \ell(h(\theta), x^\ast(\theta))$ represents the expected loss when evaluating estimator $h$ with respect to the underlying parameter distribution $D$,  and $L_T(h)=\frac{1}{m} \sum_{i=1}^m \ell(h(\theta_i), x^\ast(\theta_i))$ denotes the average empirical loss computed over the training set $T$, consisting of $m$ data samples.
The GE is a global characteristic that pertains to the class of estimators, which evaluates their suitability for given learning problem. A large generalization error indicates that within this class of estimators, there exist estimators that significantly deviate in performance on average between the true loss $L_D$, and the empirical loss $L_T$. Given the loss function $\ell$, \citet{shultzman2023generalization} proved the following theorem which gives bound for generalization error of the class of ReLU networks $\mathcal{H}^L$.
\begin{theorem}[Bound for GE, restated from \citet{shultzman2023generalization}]
    Consider the class of feed forward networks of depth-$L$ with ReLU activations 
 $\mathcal{H}^L=\{h^L(\{W^l\}_{l=1}^L) : \|W^l\|_\infty \leq B_l, l \in [1,L], \|W^1\| \leq B_1\}$ and $m$ training samples. For a loss function $\ell(h(\theta), x^\ast(\theta))=\|h(\theta) - x^\ast(\theta)\|^2$, its generalization error satisfies
    \[
    \mathcal{G}(\mathcal{H}^L) \leq 2\frac{\prod_{l=0}^L B_l}{\sqrt{m}}.
    \]
\label{thm_ge}
\end{theorem}
Thus, for $m$ training samples, it can be seen that the GE of $\mathcal{H}^L$ is bounded as $O(\frac{1}{\sqrt{m}})$ in learning our piecewise linear policy approximation.
\subsection{Approximation error.}
The AE represents the lowest possible error that an estimator can achieve within a given hypothesis class. It quantifies the amount of error incurred due to the limitation of selecting from a specific class. Unlike generalization error, the approximation error is independent of the size of the sample. For the hypothesis class $\mathcal{H}$, the AE is defined as
\[
\mathcal{A}(\mathcal{H})=\min_{h \in \mathcal{H}} L_D(h).
\]
Note that our target function $\hat{x}(\theta)$ is continuous piecewise linear function. It has been shown that any continuous piecewise linear function can be represented by a deep ReLU implementation (\citet{wang2005generalization}, \citet{arora2016understanding}), which means that the AE of NN with ReLU activations for our piecewise linear policy approximation is zero. In other words, if  $h^\ast = \arg\min_{h \in \mathcal{H}^L} L_D(h)$, $h^\ast(\theta)$ is $\varepsilon$-infeasible and $\varepsilon$-suboptimal solution.

\section{Improving Feasibility.}\label{practical}
In this section, we address issue that NN may give infeasible solutions for the original problem. 
This occurs due to the inability of NN to discern the underlying structure of an optimization problem solely through the provided training data. 
This is a critical issue because NN has been mainly used in high-assurance systems such as autonomous vehicles (\citet{kahn2017plato}), aircraft collision avoidance (\citet{julian2017neural}) and high-frequency trading (\citet{arevalo2016high}).

Therefore, we propose a strategy to improve feasibility and discuss training with suboptimal solutions.
In general, it is not easy to get an accurate optimal solution. 
Most algorithms that solve optimization problems set a certain level of threshold and stop once the threshold is achieved. 
We demonstrate that it is possible to preserve $\varepsilon$-\textit{infeasibility} even when employing suboptimal solutions for approximating parametric optimization problems, but the optimality of the approximation is inevitably reduced by the suboptimality of the suboptimal solutions.

\subsection{Infeasibility and Suboptimality.}\label{infea}
Representative approaches to improve feasibility include specifying the structure of NN, such as analyzing the output range of each layer (\citet{dutta2018output}), or obtaining a feasible solution from predictions of NN, such as greedy selection algorithm (\citet{dai2017learning}). 
In our case, we already know that an $\varepsilon$-\textit{infeasible} solution can be obtained from suitable optimization problems. 
Thus, if more strict feasible solutions than $\varepsilon$ are used, feasible solutions to the original problems can be obtained. 
Note that improving feasibility infers suboptimality. 
Therefore, before demonstrating this strategy, we discuss about feasibility and optimality of our piecewise linear policy approximation with suboptimal solutions. 
Suboptimal solution sets can be seen as a general case of optimal solution sets. If a suboptimal policy correspondence is also compact valued and upper hemicontinuous, similar arguments can be applied as in Theorem \ref{thm_main}. 
Thus, we first show this fact by starting with the following proposition.
\begin{proposition}
Define $C^\gamma(\theta)$ as $C^\gamma(\theta)=\{x\in C(\theta):|f(x,\theta)-f^\ast(\theta)|\le\gamma\}$. 
Then, $C^\gamma(\theta)$ is compact valued upper hemicontinuous correspondence.
\label{pro_sub}
\end{proposition}
\proof
We first prove Lemma \ref{lem_sub}
\begin{lemma}
If $C_1, C_2\mathrm{\,:\,}\mathrm{\Theta}\rightrightarrows X$ are correspondences, $C_1$ is upper hemicontinuous and compact valued, and $C_2$ is closed, then $C_1\cap C_2\mathrm{\,:\,}\mathrm{\Theta} \rightrightarrows X$ defined by $\ (C_1\cap C_2)(\theta)=C_1(\theta)\cap C_2(\theta)$ is upper hemicontinuous.
\label{lem_sub}
\end{lemma}
\proof[Proof of Lemma \ref{lem_sub}]
See \citet{papageorgiou1997handbook}
\endproof
To see that $C^\gamma(\theta)$ is compact valued, note that $f_\theta:C(\theta)\rightarrow\mathbb{R}$ is continuous since $f(x,\theta)$ is continuous. Since $C^\gamma(\theta)$ is closed subset of the compact set $C(\theta)$, $C^\gamma(\theta)$ is also compact. Finally, since $C(\theta)$ is compact valued continuous correspondence, it is upper hemicontinuous and compact valued. Thus, by Lemma \ref{lem_sub}, $C^\gamma(\theta)=C(\theta)\cap C^\gamma(\theta)$ is upper hemicontinuous.
\endproof
With Proposition \ref{pro_sub}, the following theorem can be proved similarly as Theorem \ref{thm_main}, which state our piecewise linear policy approximation become $\varepsilon$-infeasible and $(\varepsilon+\gamma)$-optimality solution under same conditions in Theorem \ref{thm_main}.
\begin{theorem}

Suppose that $f$, $C$ satisfy all conditions for Theorem \ref{thm_max}. 
Define a training data set $T=\{(\theta_i,x^\gamma(\theta_i))|\theta_i\in\mathrm{\Theta},i=1,\cdots,m\}$ where $x^\gamma(\theta_i)$ is an arbitrarily chosen point from $C^\gamma(\theta_i)$.
For $\theta\in Conv(\{\theta_1,\theta_2,\ldots,\theta_m\})$, define $d(S_{j_\theta})=\max\{\| p-q\|:p,q\in Conv(S_{j_\theta})\}$ where $S_{j_\theta}$ is element of finite collection $\mathcal{S}$ in Theorem \ref{thm_con}. Then, the function $\hat{x}(\theta)=\lambda^{j_\theta}_1{x^\gamma(\theta}_1^{j_\theta})+\ldots+\lambda^{j_\theta}_{k+1}{x^\gamma(\theta}_{k+1}^{j_\theta})$ satisfies the followings.
\begin{enumerate}
    \item[(a)] If $Conv(C^\gamma(\theta))\subseteq C(\theta)$, for given $\varepsilon>0$, there exists $\delta>0$ s.t. if $d(S_j)<\delta$, $\hat{x}(\theta)$ is $\varepsilon$-infeasible solution.
    \item[(b)] If $|f(x,\theta)-f^\ast(\theta)|\le\gamma, \forall x\in Conv(C^\gamma(\theta))$, for given $\varepsilon>0$, there exists $\delta>0$ s.t. if $d(S_j)<\delta$, $\hat{x}(\theta)$ is $(\varepsilon+\gamma)$-suboptimal solution
\label{thm_sub}
\end{enumerate}
\end{theorem}
\proof
(a) Since $C^\gamma$ is upper hemicontinuous, similar to statement for (a) in Theorem \ref{thm_main}, we get
\[
\lambda^{j\theta}_1{x^\gamma(\theta}_1^{j\theta})+\ldots+\lambda^{j\theta}_{k+1}{x^\gamma(\theta}_{k+1}^{j\theta})\in\mathcal{B}_\varepsilon(Conv(C^\gamma(\theta))).
\]
With the assumption, $\hat{x}(\theta)\in\mathcal{B}_\varepsilon(Conv(C^\gamma(\theta)))\subseteq\mathcal{B}_\varepsilon(C(\theta))$.

(b) Assume $|f(x,\theta)-f^\ast(\theta)|\le\gamma, \forall x\in Conv(C^\gamma(\theta))$ and let $\varepsilon>0$. Since $C^\gamma$ is upper hemicontinuous and compact valued, similar to statement for (b) in Theorem \ref{thm_main}, there exists $\delta^\prime>0$ such that, 
\[
\inf_{x^\gamma\in Conv(C^\gamma(\theta))} \|x^\prime-x^\gamma\|<\delta^\prime\Longrightarrow |f(x^\prime,\theta)-f_{min}^\gamma(\theta)|<\varepsilon, \forall x^\prime\in\mathcal{B}_\varepsilon(C(\theta)).
\]
where $x_{min}^\gamma={\arg\min}_{x^\gamma\in Conv(C^\gamma(\theta))}\|x^\prime-x^\gamma\|$. 
Since $x_{min}^\gamma\in Conv(C^\gamma(\theta))$, we have
\[
|f(x^\prime,\theta)-f^\ast(\theta)|\leq|f(x^\prime,\theta)-f(x^\gamma_{min},\theta)|+|f(x^\gamma_{min},\theta)-f^\ast(\theta)|<\varepsilon+\gamma.
\]
Then, similar to (b) of Theorem \ref{thm_main}, we have
\[
\hat{x}(\theta)=\lambda^{j\theta}_1{x^\gamma(\theta}_1^{j\theta})+\ldots+\lambda^{j\theta}_{k+1}{x^\gamma(\theta}_{k+1}^{j\theta})\in\mathcal{B}_{\delta^\prime}(Conv(C^\gamma(\theta))).
\]
Also, note that condition in (b) implies $Conv(C^\gamma(\theta))\subseteq C^\gamma(\theta)\subseteq C(\theta)$ which indicates $\hat{x}(\theta)\in\mathcal{B}_\varepsilon(C(\theta))$. Accordingly, $|f(\hat{x}(\theta),\theta)-f^\ast(\theta)|<\varepsilon+\gamma$.

\endproof

Now we demonstrate the proposed strategy with linear programming (LP) and quadratic programming (QP). Note that, since the optimal solution set of every convex optimization is convex, convex optimization satisfies the condition for part (a) in Theorem \ref{thm_main}. Thus, our strategy for improving feasibility can be applied. The formulations of LP and QP are as follows which are modified version with with slightly perturbed right hand side in standard formulation. For each problem, parameters except $\theta$ and $t$ were randomly generated and fixed.

For LP,
\[
\min_x \quad c^{T}x \quad \text{subject\;to}\quad Ax \leq \theta-t\mathbf{1},\quad c, x,\theta\in\mathbb{R}^n, \quad A\in\mathbb{R}^{nxn} \tag{2} \label{eq:2}
\]
For QP,
\[
\min_x \quad (1/2)x^{T}Px+q^{T}x \quad \text{subject\;to}\quad Gx \leq \theta-t\mathbf{1},\quad q, x,\theta\in\mathbb{R}^n, \quad P, G \in \mathbb{R}^{nxn} \tag{3} \label{eq:3}
\] 
where $\mathbf{1}=
\begin{bmatrix}
1\\
1\\
\vdots\\
1
\end{bmatrix}$

\begin{figure}
     \centering
     \begin{subfigure}[t]{0.96\textwidth}
         \centering
         \includegraphics[width=\textwidth]{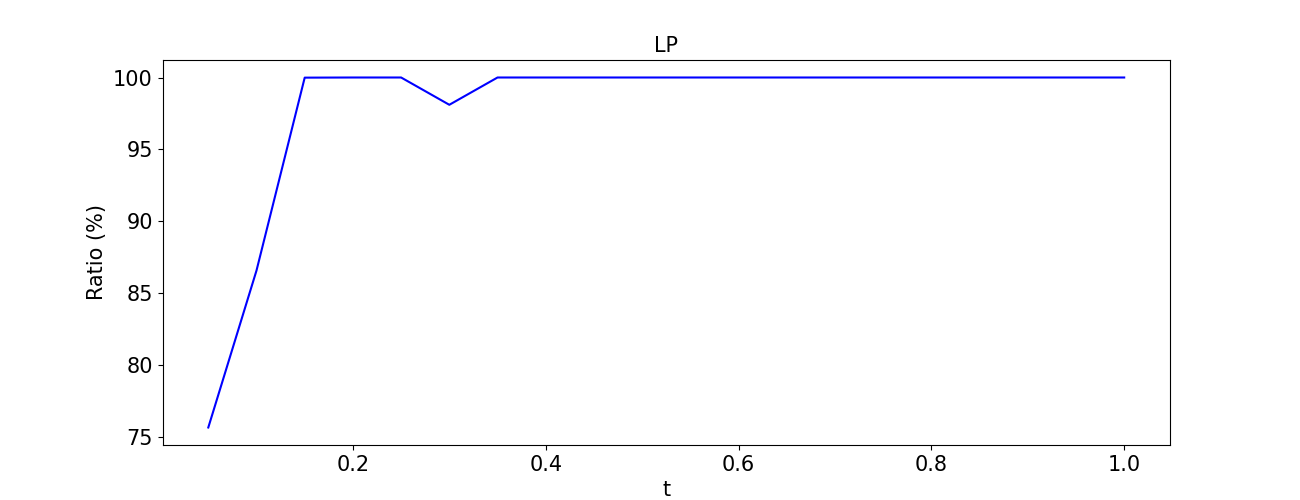}
     \end{subfigure}
     \hfill
     \begin{subfigure}[t]{0.96\textwidth}
         \centering
         \includegraphics[width=\textwidth]{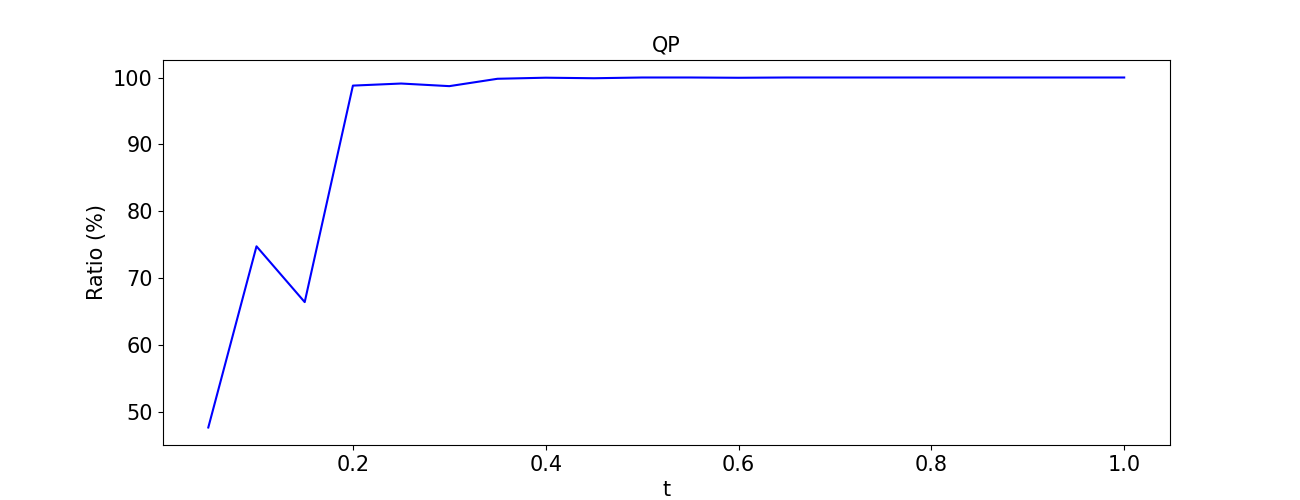}
     \end{subfigure}
        \caption{Ratio(\%) of feasible solutions to original problem ($t=0$)}
        \label{fig:4}
\end{figure}

Here, $t$ serves to obtain a solution further inside than the feasible region of the original problem. Problems (\ref{eq:2}) and (\ref{eq:3}) are solved for a total of 10,000 times each while slightly increasing the value of $t$ for each iteration. At this time, 10,000 samples for $\theta$ are generated from standard normal distribution for each problem. For each $t$, we trained NN with training pair $\{\theta,x(\theta)\}$ and calculated the ratio of the feasible approximation of NN to the problem when $t=0$. As shown in Figure \ref{fig:4}, the ratio converges to 100\% in every problem, which indicates our strategy guarantees feasibility.

\section{Conclusion.}\label{conclu}
In this paper, we build theoretical foundations on approximating direct mapping from
input parameters to optimal solution through NN universally and derive conditions that allow the Universal Approximation Theorem to be applied to parametric optimization problems by constructing piecewise linear policy approximation explicitly. More specifically, we cast single-valued continuous piecewise linear approximation for optimal solution of parametric optimization and analyze it in terms of feasibility and optimality and show that NN with ReLU activations can be valid approximator in terms of generalization and approximation error. Moreover, we propose strategy to improve feasibility and discuss on the suboptimal training data, findings from this study can directly benefit solving parametric optimization problems in real-time control systems or high-assurance systems. In future research, we plan to extend our theory to more general parametric optimization problems such as integer programming, and study more approaches for addressing infeasibility of approximated solutions.

\bibliographystyle{unsrtnat}
\bibliography{UAPO}  






\end{document}